# CONSISTENT ESTIMATION OF THE BASIC NEIGHBORHOOD OF MARKOV RANDOM FIELDS

By Imre Csiszár[1] and Zsolt Talata [2]

*Hungarian Academy of Sciences*

For Markov random fields on $\mathbb{Z}^d$ with finite state space, we address the statistical estimation of the basic neighborhood, the smallest region that determines the conditional distribution at a site on the condition that the values at all other sites are given. A modification of the Bayesian Information Criterion, replacing likelihood by pseudo-likelihood, is proved to provide strongly consistent estimation from observing a realization of the field on increasing finite regions: the estimated basic neighborhood equals the true one eventually almost surely, not assuming any prior bound on the size of the latter. Stationarity of the Markov field is not required, and phase transition does not affect the results.

**1. Introduction.** In this paper Markov random fields on the lattice $\mathbb{Z}^d$ with finite state space are considered, adopting the usual assumption that the finite-dimensional distributions are strictly positive. Equivalently, these are Gibbs fields with finite range interaction; see [13]. They are essential in statistical physics, for modeling interactive particle systems [10], and also in several other fields [3], for example, in image processing [2].

One statistical problem for Markov random fields is parameter estimation when the interaction structure is known. By this we mean knowledge of the *basic neighborhood*, the minimal lattice region that determines the conditional distribution at a site on the condition that the values at all other sites are given; formal definitions are in Section 2. The conditional probabilities involved, assumed translation invariant, are parameters of the model. Note

[1]Supported by Hungarian National Foundation for Scientific Research Grants T26041, T32323, TS40719 and T046376.

[2]Supported by Hungarian National Foundation for Scientific Research Grant T046376.

that they need not uniquely determine the joint distribution on $\mathbb{Z}^d$, a phenomenon known as *phase transition*. Another statistical problem is *model selection*, that is, the statistical estimation of the interaction structure (the basic neighborhood). This paper is primarily devoted to the latter.

Parameter estimation for Markov random fields with a known interaction structure was considered by, among others, Pickard [19], Gidas [14, 15], Geman and Graffigne [12] and Comets [6]. Typically, parameter estimation does not directly address the conditional probabilities mentioned above, but rather the *potential*. This admits parsimonious representation of the conditional probabilities that are not free parameters, but have to satisfy algebraic conditions that need not concern us here. For our purposes, however, potentials will not be needed.

We are not aware of papers addressing model selection in the context of Markov random fields. In other contexts, penalized likelihood methods are popular; see [1, 21]. The Bayesian Information Criterion (BIC) of Schwarz [21] has been proven to lead to consistent estimation of the "order of the model" in various cases, such as i.i.d. processes with distributions from exponential families [17], autoregressive processes [16] and Markov chains [11]. These proofs include the assumption that the number of candidate model classes is finite; for Markov chains this means that there is a known upper bound on the order of the process. The consistency of the BIC estimator of the order of a Markov chain without such prior bound was proved by Csiszár and Shields [8]; further related results appear in [7]. A related recent result, for processes with variable memory length [5, 22], is the consistency of the BIC estimator of the context tree, without any prior bound on memory depth [9].

For Markov random fields, penalized likelihood estimators like BIC run into the problem that the likelihood function cannot be calculated explicitly. In addition, no simple formula is available for the "number of free parameters" typically used in the penalty term. To overcome these problems, we will replace likelihood by pseudo-likelihood, first introduced by Besag [4], and modify also the penalty term; this will lead us to an analogue of BIC called the *Pseudo-Bayesian Information Criterion* or PIC. Our main result is that if one minimizes this criterion for a family of hypothetical basic neighborhoods that grows with the sample size at a specified rate, the resulting PIC estimate of the basic neighborhood equals the true one eventually almost surely. In particular, the consistency theorem does not require a prior upper bound on the size of the basic neighborhood. It should be emphasized that the underlying Markov field need not be stationary (translation invariant), and phase transition causes no difficulty.

An auxiliary result perhaps of independent interest is a typicality proposition on the uniform closeness of empirical conditional probabilities to the true ones, for conditioning regions whose size may grow with the sample size. Though this result is weaker than analogous ones for Markov chains in [7], it will be sufficient for our purposes.

The structure of the paper is the following. In Section 2 we introduce the basic notation and definitions, and formulate the main result. Its proof is provided by the propositions in Sections 4 and 5. Section 3 contains the statement and proof of the typicality proposition. Section 4 excludes overestimation, that is, the possibility that the estimated basic neighborhood properly contains the true one, using the typicality proposition. Section 5 excludes underestimation, that is, the possibility that the estimated basic neighborhood does not contain the true one, via an entropy argument and a modification of the typicality result. Section 6 is a discussion of the results. The Appendix contains some technical lemmas.

**2. Notation and statement of the main results.** We consider the $d$-dimensional *lattice* $\mathbb{Z}^d$. The points $i \in \mathbb{Z}^d$ are called sites, and $\|i\|$ denotes the maximum norm of $i$, that is, the maximum of the absolute values of the coordinates of $i$. The cardinality of a finite set $\Delta$ is denoted by $|\Delta|$. The notation $\subseteq$ and $\subset$ of inclusion and strict inclusion are distinguished in this paper.

A *random field* is a family of random variables indexed by the sites of the lattice, $\{X(i) : i \in \mathbb{Z}^d\}$, where each $X(i)$ is a random variable with values in a finite set $A$. For $\Delta \subseteq \mathbb{Z}^d$, a region of the lattice, we write $X(\Delta) = \{X(i) : i \in \Delta\}$. For the realizations of $X(\Delta)$ we use the notation $a(\Delta) = \{a(i) \in A : i \in \Delta\}$. When $\Delta$ is finite, the $|\Delta|$-tuples $a(\Delta) \in A^\Delta$ will be referred to as *blocks*.

The joint distribution of the random variables $X(i)$ is denoted by $Q$. We assume that its finite-dimensional marginals are strictly positive, that is,

$$Q(a(\Delta)) = \text{Prob}\{X(\Delta) = a(\Delta)\} > 0 \qquad \text{for } \Delta \subset \mathbb{Z}^d \text{ finite}, a(\Delta) \in A^\Delta.$$

The last standard assumption admits unambiguous definition of the conditional probabilities

$$Q(a(\Delta)|a(\Phi)) = \text{Prob}\{X(\Delta) = a(\Delta)|X(\Phi) = a(\Phi)\}$$

for all disjoint finite regions $\Delta$ and $\Phi$.

By a *neighborhood* $\Gamma$ (of the origin 0) we mean a finite, central-symmetric set of sites with $0 \notin \Gamma$. Its radius is $r(\Gamma) = \max_{i\in\Gamma} \|i\|$. For any $\Delta \subseteq \mathbb{Z}^d$, its translate when 0 is translated to $i$ is denoted by $\Delta^i$. The translate $\Gamma^i$ of a neighborhood $\Gamma$ (of the origin) will be called the $\Gamma$-neighborhood of the site $i$; see Figure 1.

A *Markov random field* is a random field as above such that there exists a neighborhood $\Gamma$, called a *Markov neighborhood*, satisfying for every $i \in \mathbb{Z}^d$

$$Q(a(i)|a(\Delta^i)) = Q(a(i)|a(\Gamma^i)) \qquad \text{if } \Delta \supset \Gamma, 0 \notin \Delta, \tag{2.1}$$

where the last conditional probability is translation invariant.

This concept is equivalent to that of a Gibbs field with a finite range interaction; see [13]. Motivated by this fact, the matrix

$$Q_\Gamma = \{Q_\Gamma(a|a(\Gamma)) : a \in A, a(\Gamma) \in A^\Gamma\}$$

specifying the (positive, translation-invariant) conditional probabilities in (2.1) will be called *one-point specification*. All distributions on $A^{\mathbb{Z}^d}$ that satisfy (2.1) with a given conditional probability matrix $Q_\Gamma$ are called *Gibbs distributions* with one-point specification $Q_\Gamma$. The distribution $Q$ of the given Markov random field is one of these; $Q$ is not necessarily translation invariant.

The following lemma summarizes some well-known facts; their formal derivation from results in [13] is indicated in the Appendix.

Lemma 2.1. *For a Markov random field on the lattice as above, there exists a neighborhood $\Gamma_0$ such that the Markov neighborhoods are exactly those that contain $\Gamma_0$. Moreover, the global Markov property*

$$Q(a(\Delta)|a(\mathbb{Z}^d \setminus \Delta)) = Q\left(a(\Delta)\middle|a\left(\bigcup_{i\in\Delta} \Gamma_0^i \setminus \Delta\right)\right)$$

*holds for each finite region $\Delta \subset \mathbb{Z}^d$. These conditional probabilities are translation invariant and uniquely determined by the one-point specification $Q_{\Gamma_0}$.*

The smallest Markov neighborhood $\Gamma_0$ of Lemma 2.1 will be called the *basic neighborhood*. The minimal element of the corresponding one-point specification matrix $Q_{\Gamma_0}$ is denoted by $q_{\min}$:

$$q_{\min} = \min_{a\in A, a(\Gamma_0)\in A^{\Gamma_0}} Q_{\Gamma_0}(a|a(\Gamma_0)) > 0.$$

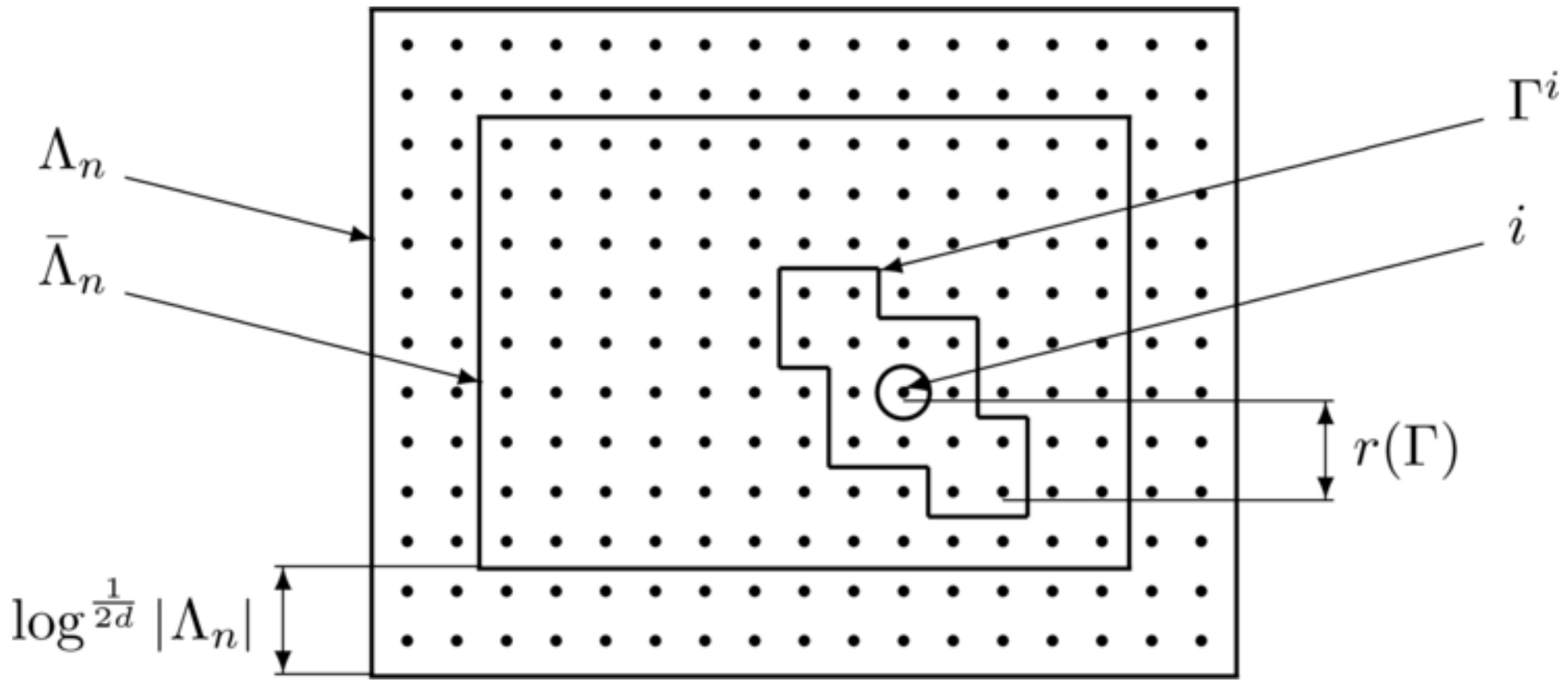

Fig. 1. *The $\Gamma$-neighborhood of the site $i$, and the sample region $\Lambda_n$.*

In this paper we are concerned with the statistical estimation of the basic neighborhood $\Gamma_0$ from observation of a realization of the Markov random field on an increasing sequence of finite regions $\Lambda_n \subset \mathbb{Z}^d$, $n \in \mathbb{N}$; thus the $n$th sample is $x(\Lambda_n)$.

We will draw the statistical inference about a possible basic neighborhood $\Gamma$ based on the blocks $a(\Gamma) \in A^\Gamma$ appearing in the sample $x(\Lambda_n)$. For technical reasons, we will consider only such blocks whose center is in a subregion $\bar{\Lambda}_n$ of $\Lambda_n$, consisting of those sites $i \in \Lambda_n$ for which the ball with center $i$ and radius $\log^{1/(2d)} |\Lambda_n|$ also belongs to $\Lambda_n$:

$$\bar{\Lambda}_n = \{i \in \Lambda_n : \{j \in \mathbb{Z}^d : \|i - j\| \le \log^{1/(2d)} |\Lambda_n|\} \subseteq \Lambda_n\};$$

see Figure 1. Our only assumptions about the sample regions $\Lambda_n$ will be that

$$\Lambda_1 \subset \Lambda_2 \subset \cdots; \qquad |\Lambda_n|/|\bar{\Lambda}_n| \to 1.$$

For each block $a(\Gamma) \in A^\Gamma$, let $N_n(a(\Gamma))$ denote the number of occurrences of the block $a(\Gamma)$ in the sample $x(\Lambda_n)$ with the center in $\bar{\Lambda}_n$,

$$N_n(a(\Gamma)) = |\{i \in \bar{\Lambda}_n : \Gamma^i \subseteq \Lambda_n, x(\Gamma^i) = a(\Gamma)\}|.$$

The blocks corresponding to $\Gamma$-neighborhoods completed with their centers will be denoted briefly by $a(\Gamma, 0)$. Similarly as above, for each $a(\Gamma, 0) \in A^{\Gamma \cup \{0\}}$ we write

$$N_n(a(\Gamma, 0)) = |\{i \in \bar{\Lambda}_n : \Gamma^i \subseteq \Lambda_n, x(\Gamma^i \cup \{i\}) = a(\Gamma, 0)\}|.$$

The notation $a(\Gamma, 0) \in x(\Lambda_n)$ will mean that $N_n(a(\Gamma, 0)) \ge 1$.

The restriction $\Gamma^i \subseteq \Lambda_n$ in the above definitions is automatically satisfied if $r(\Gamma) \le \log^{1/(2d)} |\Lambda_n|$. Hence the same number of blocks is taken into account for all neighborhoods, except for very large ones:

$$\sum_{a(\Gamma) \in A^\Gamma} N_n(a(\Gamma)) = |\bar{\Lambda}_n| \qquad \text{if } r(\Gamma) \le \log^{1/(2d)} |\Lambda_n|.$$

For Markov random fields the likelihood function cannot be explicitly determined. We shall use instead the pseudo-likelihood defined below.

Given the sample $x(\Lambda_n)$, the *pseudo-likelihood* function associated with a neighborhood $\Gamma$ is the following function of a matrix $Q'_\Gamma$ regarded as the one-point specification of a hypothetical Markov random field for which $\Gamma$ is a Markov neighborhood:

$$\begin{aligned} \mathrm{PL}_\Gamma(x(\Lambda_n), Q'_\Gamma) &= \prod_{i \in \bar{\Lambda}_n} Q'_\Gamma(x(i) | x(\Gamma^i)) \\ &= \prod_{a(\Gamma,0) \in x(\Lambda_n)} Q'_\Gamma(a(0) | a(\Gamma))^{N_n(a(\Gamma,0))}. \end{aligned} \tag{2.2}$$

We note that not all matrices $Q'_\Gamma$ satisfying

$$\sum_{a\in A} Q'_\Gamma(a(0)|a(\Gamma)) = 1, \qquad a(\Gamma) \in A^\Gamma$$

are possible one-point specifications; the elements of a one-point specification matrix have to satisfy several algebraic relations not shown here. Still, we define the pseudo-likelihood also for $Q'_\Gamma$ not satisfying those relations, even admitting some elements of $Q'_\Gamma$ to be 0.

The maximum of this pseudo-likelihood is attained for $Q'_\Gamma(a(0)|a(\Gamma)) = \frac{N_n(a(\Gamma,0))}{N_n(a(\Gamma))}$. Thus, given the sample $x(\Lambda_n)$, the logarithm of the *maximum pseudo-likelihood* for the neighborhood $\Gamma$ is

$$\log \mathrm{MPL}_\Gamma(x(\Lambda_n)) = \sum_{a(\Gamma,0)\in x(\Lambda_n)} N_n(a(\Gamma,0)) \log \frac{N_n(a(\Gamma,0))}{N_n(a(\Gamma))}. \tag{2.3}$$

Now we are able to formalize a criterion in analogy to the Bayesian Information Criterion that can be calculated from the sample.

Definition 2.1. Given a sample $x(\Lambda_n)$, the Pseudo-Bayesian Information Criterion, in short PIC, for the neighborhood $\Gamma$ is

$$\mathrm{PIC}_\Gamma(x(\Lambda_n)) = -\log \mathrm{MPL}_\Gamma(x(\Lambda_n)) + |A|^{|\Gamma|} \log |\Lambda_n|.$$

Remark. In our penalty term, the number $|A|^{|\Gamma|}$ of possible blocks $a(\Gamma) \in A^\Gamma$ replaces "half the number of free parameters" appearing in BIC, for which number no simple formula is available. Note that our results remain valid, with the same proofs, if the above penalty term is multiplied by any $c > 0$.

The PIC estimator of the basic neighborhood $\Gamma_0$ is defined as that hypothetical $\Gamma$ for which the value of the criterion is minimal. An important feature of our estimator is that the family of hypothetical $\Gamma$'s is allowed to extend as $n \to \infty$, and thus no *a priori* upper bound for the size of the unknown $\Gamma_0$ is needed. Our main result says the PIC estimator is strongly consistent if the hypothetical $\Gamma$'s are those with $r(\Gamma) \le r_n$, where $r_n$ grows sufficiently slowly.

We mean by strong consistency that the estimated basic neighborhood equals $\Gamma_0$ eventually almost surely as $n \to \infty$. Here and in the sequel, "eventually almost surely" means that with probability 1 there exists a threshold $n_0$ [depending on the infinite realization $x(\mathbb{Z}^d)$] such that the claim holds for all $n \ge n_0$.

THEOREM 2.1. *The PIC estimator*

$$\widehat{\Gamma}_{\mathrm{PIC}}(x(\Lambda_n)) = \mathop{\arg\min}_{\Gamma:\, r(\Gamma)\le r_n} \mathrm{PIC}_\Gamma(x(\Lambda_n)),$$

*with*

$$r_n = o(\log^{1/(2d)}|\Lambda_n|),$$

*satisfies*

$$\widehat{\Gamma}_{\mathrm{PIC}}(x(\Lambda_n)) = \Gamma_0$$

*eventually almost surely as* $n \to \infty$.

PROOF. Theorem 2.1 follows from Propositions 4.1 and 5.1 below. □

REMARK. Actually, the assertion will be proved for $r_n$ equal to a constant times $\log^{1/(2d)}|\bar{\Lambda}_n|$. However, as this constant depends on the unknown distribution $Q$, the consistency can be guaranteed only when

$$r_n = o(\log^{1/(2d)}|\bar{\Lambda}_n|) = o(\log^{1/(2d)}|\Lambda_n|).$$

It remains open whether consistency holds when the hypothetical neighborhoods are allowed to grow faster, or even without any condition on the hypothetical neighborhoods.

As a consequence of the above, we are able to construct a strongly consistent estimator of the one-point specification $Q_{\Gamma_0}$.

COROLLARY 2.1. *The empirical estimator of the one-point specification,*

$$\widehat{Q}_{\widehat{\Gamma}}(a(0)|a(\widehat{\Gamma})) = \frac{N_n(a(\widehat{\Gamma},0))}{N_n(a(\widehat{\Gamma}))}, \qquad a(0)\in A, a(\widehat{\Gamma})\in A^{\widehat{\Gamma}},$$

*converges to the true* $Q_{\Gamma_0}$ *almost surely as* $n \to \infty$, *where* $\widehat{\Gamma}$ *is the PIC estimator* $\widehat{\Gamma}_{\mathrm{PIC}}$.

PROOF. Immediate from Theorem 2.1 and Proposition 3.1 below. □

## 3. The typicality result.

PROPOSITION 3.1. *Simultaneously for all Markov neighborhoods with* $r(\Gamma) \le \alpha^{1/(2d)}\log^{1/(2d)}|\bar{\Lambda}_n|$ *and blocks* $a(\Gamma,0)\in A^{\Gamma\cup\{0\}}$,

$$\left|\frac{N_n(a(\Gamma,0))}{N_n(a(\Gamma))} - Q(a(0)|a(\Gamma))\right| < \sqrt{\frac{\kappa\log N_n(a(\Gamma))}{N_n(a(\Gamma))}}$$

*eventually almost surely as* $n\to\infty$, *if*

$$0<\alpha\le 1, \qquad \kappa > 2^{3d}e\alpha\log(|A|^2+1).$$

To prove this proposition we will use an idea similar to the "coding technique" of Besag [3]; namely, we partition $\bar{\Lambda}_n$ into subsets $\bar{\Lambda}_n^k$ such that the random variables at the sites $i \in \bar{\Lambda}_n^k$ are conditionally independent given the values of those at the other sites. First we introduce some further notation. Let

$$R_n = \lfloor \alpha^{1/(2d)} \lceil \log |\bar{\Lambda}_n| \rceil^{1/(2d)} \rfloor. \tag{3.4}$$

We partition the region $\bar{\Lambda}_n$ by intersecting it with sublattices of $\mathbb{Z}^d$ such that the distance between sites in a sublattice is $4R_n + 1$. The intersections of $\bar{\Lambda}_n$ with these sublattices will be called sieves. Indexed by the offset $k$ relative to the origin 0, the sieves are

$$\bar{\Lambda}_n^k = \{i \in \bar{\Lambda}_n : i = k + (4R_n + 1)v, v \in \mathbb{Z}^d\}, \qquad \|k\| \leq 2R_n;$$

see Figure 2. For a neighborhood $\Gamma$, let $N_n^k(a(\Gamma))$ denote the number of occurrences of the block $a(\Gamma) \in A^\Gamma$ in the sample $x(\Lambda_n)$ with center in $\bar{\Lambda}_n^k$,

$$N_n^k(a(\Gamma)) = |\{i \in \bar{\Lambda}_n^k : \Gamma^i \subseteq \Lambda_n, x(\Gamma^i) = a(\Gamma)\}|.$$

Similarly, let

$$N_n^k(a(\Gamma, 0)) = |\{i \in \bar{\Lambda}_n^k : \Gamma^i \subseteq \Lambda_n, x(\Gamma^i \cup \{i\}) = a(\Gamma, 0)\}|.$$

Clearly,

$$N_n(a(\Gamma)) = \sum_{k : \|k\| \leq 2R_n} N_n^k(a(\Gamma)) \quad \text{and} \quad N_n(a(\Gamma, 0)) = \sum_{k : \|k\| \leq 2R_n} N_n^k(a(\Gamma, 0)).$$

The notation $a(\Gamma) \in x(\Lambda_n^k)$ will mean that $N_n^k(a(\Gamma)) \geq 1$.

Denote by $\Phi_n(\Gamma)$ the set of sites outside the neighborhood $\Gamma$ whose norm is at most $2R_n$,

$$\Phi_n(\Gamma) = \{i \in \mathbb{Z}^d : \|i\| \leq 2R_n, i \notin \Gamma\};$$

see Figure 2. $\Phi_n^i(\Gamma)$ denotes the translate of $\Phi_n(\Gamma)$ when 0 is translated to $i$.

For a finite region $\Xi \subset \mathbb{Z}^d$, conditional probabilities on the condition $X(\Xi) = x(\Xi) \in A^\Xi$ will be denoted briefly by $\text{Prob}\{\cdot \mid x(\Xi)\}$.

In the following lemma the neighborhoods $\Gamma$ need not be Markov neighborhoods.

Lemma 3.1. *Simultaneously for all sieves $k$, neighborhoods $\Gamma$ with $r(\Gamma) \leq R_n$ and blocks $a(\Gamma) \in A^\Gamma$,*

$$(1 + \varepsilon) \log N_n^k(a(\Gamma)) \geq \log |\bar{\Lambda}_n|,$$

*eventually almost surely as $n \to \infty$, where $\varepsilon > 0$ is an arbitrary constant.*

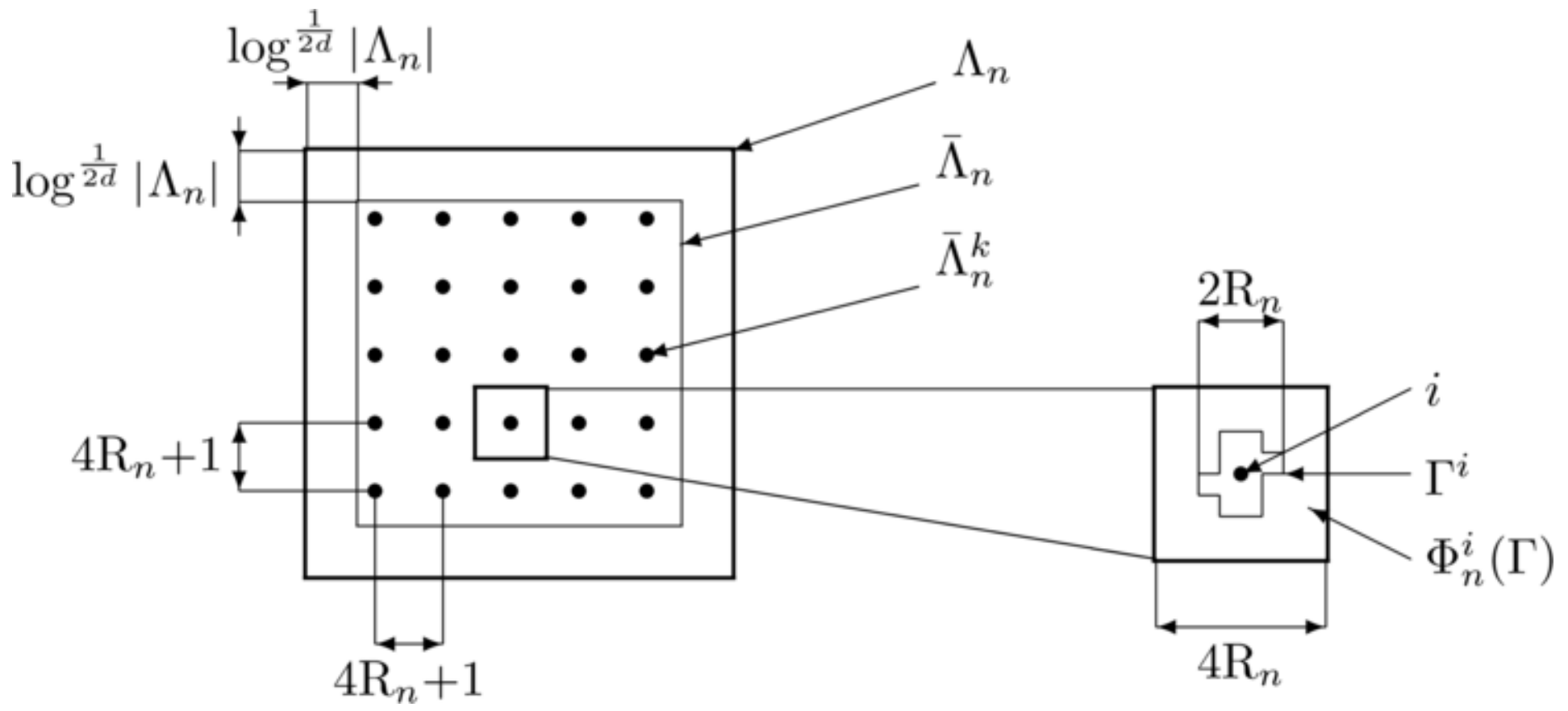

Fig. 2. *The sieve* $\bar{\Lambda}_n^k$.

Proof. As a consequence of Lemma 2.1, for any fixed sieve $k$ and neighborhood $\Gamma$ with $r(\Gamma) \le R_n$, the random variables $X(\Gamma^i)$, $i \in \bar{\Lambda}_n^k$, are conditionally independent given the values of the random variables in the rest of the sites of the sample region $\Lambda_n$. By Lemma A.5 in the Appendix,

$$Q(a(\Gamma)|a(\Phi_n(\Gamma))) \ge q_{\min}^{|\Gamma|}, \qquad a(\Phi_n(\Gamma)) \in A^{\Phi_n(\Gamma)},$$

hence we can use the large deviation theorem of Lemma A.3 in the Appendix with $p_* = q_{\min}^{|\Gamma|}$ to obtain

$$\mathrm{Prob}\left\{\frac{N_n^k(a(\Gamma))}{|\bar{\Lambda}_n^k|} < \frac{1}{2}q_{\min}^{|\Gamma|}\middle| x\left(\Lambda_n \Big\backslash \bigcup_{i\in\bar{\Lambda}_n^k} \Gamma^i\right)\right\} \le \exp\left[-|\bar{\Lambda}_n^k|\frac{q_{\min}^{|\Gamma|}}{16}\right].$$

Hence also for the unconditional probabilities,

$$\mathrm{Prob}\left\{\frac{N_n^k(a(\Gamma))}{|\bar{\Lambda}_n^k|} < \frac{1}{2}q_{\min}^{|\Gamma|}\right\} \le \exp\left[-|\bar{\Lambda}_n^k|\frac{q_{\min}^{|\Gamma|}}{16}\right].$$

Note that for $n \ge n_0$ (not depending on $k$) we have

$$|\bar{\Lambda}_n^k| \ge \frac{1}{2}\frac{|\bar{\Lambda}_n|}{(4R_n+1)^d} > \frac{|\bar{\Lambda}_n|}{(5R_n)^d}.$$

Using this and the consequence $|\Gamma| \le (2R_n+1)^d < (3R_n)^d$ of $r(\Gamma) \le R_n$, the last probability bound implies for $n \ge n_0$

$$\mathrm{Prob}\left\{\frac{N_n^k(a(\Gamma))}{|\bar{\Lambda}_n|} < \frac{q_{\min}^{(3R_n)^d}}{2(5R_n)^d}\right\} \le \exp\left[-|\bar{\Lambda}_n|\frac{q_{\min}^{(3R_n)^d}}{16(5R_n)^d}\right].$$

Using the union bound and Lemma A.6 in the Appendix, it follows that

$$\begin{aligned}\mathrm{Prob}\Big\{&\frac{N_n^k(a(\Gamma))}{|\bar\Lambda_n|} < \frac{q_{\min}^{(3R_n)^d}}{2(5R_n)^d},\\ &\text{for some } k,\Gamma,a(\Gamma) \text{ with } \|k\|\le 2R_n, r(\Gamma)\le R_n, a(\Gamma)\in A^\Gamma\Big\}\\ &\le \exp\left[-|\bar\Lambda_n|\frac{q_{\min}^{(3R_n)^d}}{16(5R_n)^d}\right]\cdot(4R_n+1)^d\cdot(|A|^2+1)^{(2R_n+1)^d/2}.\end{aligned}$$

Recalling (3.4), this is summable in $n$, and thus the Borel–Cantelli lemma gives

$$N_n^k(a(\Gamma)) \ge |\bar\Lambda_n|\frac{q_{\min}^{3^d\alpha^{1/2}(1+\log|\bar\Lambda_n|)^{1/2}}}{2\cdot 5^d\alpha^{1/2}(1+\log|\bar\Lambda_n|)^{1/2}},$$

eventually almost surely as $n\to\infty$, simultaneously for all sieves $k$, neighborhoods $\Gamma$ with $r(\Gamma)\le R_n$ and blocks $a(\Gamma)\in A^\Gamma$. This proves the lemma. □

Lemma 3.2. *Simultaneously for all sieves $k$, Markov neighborhoods $\Gamma$ with $r(\Gamma)\le R_n$ and blocks $a(\Gamma,0)\in A^{\Gamma\cup\{0\}}$,*

$$\left|\frac{N_n^k(a(\Gamma,0))}{N_n^k(a(\Gamma))} - Q(a(0)|a(\Gamma))\right| < \sqrt{\frac{\delta\log^{1/2}N_n^k(a(\Gamma))}{N_n^k(a(\Gamma))}},$$

*eventually almost surely as $n\to\infty$, if*

$$\delta > 2^d e\alpha^{1/2}\log(|A|^2+1).$$

Proof. Given a sieve $k$, a Markov neighborhood $\Gamma$ and a block $a(\Gamma,0)$, the difference $N_n^k(a(\Gamma,0)) - N_n^k(a(\Gamma))Q(a(0)|a(\Gamma))$ equals

$$Y_n = \sum_{i\in\bar\Lambda_n^k\,:\,x(\Gamma^i)=a(\Gamma)} [\mathbb{I}(X(i)=a(0)) - Q(a(0)|a(\Gamma))],$$

where $\mathbb{I}(\cdot)$ denotes the indicator function; hence the claimed inequality is equivalent to

$$-\sqrt{N_n^k(a(\Gamma))\delta\log^{1/2}N_n^k(a(\Gamma))} < Y_n < \sqrt{N_n^k(a(\Gamma))\delta\log^{1/2}N_n^k(a(\Gamma))}.$$

We will prove that the last inequalities hold eventually almost surely as $n\to\infty$, simultaneously for all sieves $k$, Markov neighborhoods $\Gamma$ with $r(\Gamma)\le R_n$ and blocks $a(\Gamma,0)\in A^{\Gamma\cup\{0\}}$. We concentrate on the second inequality; the proof for the first one is similar.

Denote

$$G_j(k,a(\Gamma,0)) = \left\{ \max_{n\in\mathcal{N}_j(k,a(\Gamma))} Y_n \geq \sqrt{e^j \delta j^{1/2}} \right\},$$

where

$$\mathcal{N}_j(k,a(\Gamma)) = \{n : e^j < N_n^k(a(\Gamma)) \leq e^{j+1}, (1+\varepsilon)\log N_n^k(a(\Gamma)) \geq \log|\bar{\Lambda}_n|\};$$

if $n\in\mathcal{N}_j(k,a(\Gamma))$, then by (3.4)

$$R_n = \lfloor \alpha^{1/(2d)} \lceil \log|\bar{\Lambda}_n| \rceil^{1/(2d)} \rfloor \leq \alpha^{1/(2d)}(1+(1+\varepsilon)(j+1))^{1/(2d)} \stackrel{\text{def}}{=} R^{(j)}. \tag{3.5}$$

The claimed inequality $Y_n < \sqrt{N_n^k(a(\Gamma))\delta \log^{1/2} N_n^k(a(\Gamma))}$ holds for each $n$ with $e^j < N_n^k(a(\Gamma)) \leq e^{j+1}$ if

$$\max_{n : e^j < N_n^k(a(\Gamma)) \leq e^{j+1}} Y_n < \sqrt{e^j \delta j^{1/2}}.$$

By Lemma 3.1, the condition $(1+\varepsilon)\log N_n^k(a(\Gamma)) \geq \log|\bar{\Lambda}_n|$ in the definition of $\mathcal{N}_j(k,a(\Gamma))$ is satisfied eventually almost surely, simultaneously for all sieves $k$, neighborhoods $\Gamma$ with $r(\Gamma) \leq R_n$ and blocks $a(\Gamma) \in A^\Gamma$. Hence it suffices to prove that the following holds with probability 1: the union of the events $G_j(k,a(\Gamma,0))$ for all $k$ with $\|k\| \leq 2R^{(j)}$, all $\Gamma \supseteq \Gamma_0$ with $r(\Gamma) \leq R^{(j)}$ and all $a(\Gamma,0) \in A^{\Gamma\cup\{0\}}$, obtains only for finitely many $j$.

As $n\in\mathcal{N}_j(k,a(\Gamma))$ implies $j < \log|\bar{\Lambda}_n| \leq (1+\varepsilon)(j+1)$,

$$G_j(k,a(\Gamma,0)) \subseteq \bigcup_{l=j}^{\lfloor (1+\varepsilon)(j+1) \rfloor} \left\{ \max_{n\in\mathcal{N}_{j,l}(k,a(\Gamma))} Y_n \geq \sqrt{e^j \delta j^{1/2}} \right\}, \tag{3.6}$$

where

$$\mathcal{N}_{j,l}(k,a(\Gamma)) = \{n : e^j < N_n^k(a(\Gamma)) \leq e^{j+1}, l < \log|\bar{\Lambda}_n| \leq l+1\}.$$

The random variables $X(i)$, $i\in\bar{\Lambda}_n^k$, are conditionally independent given the values of the random variables in their $\Gamma$-neighborhoods. Moreover, those $X(i)$'s for which the same block $a(\Gamma)$ appears in their $\Gamma$-neighborhood are also conditionally i.i.d. Hence $Y_n$ is the sum of $N_n^k(a(\Gamma))$ conditionally i.i.d. random variables with mean 0 and variance

$$\tfrac{1}{4} \geq D^2 = Q(a(0)|a(\Gamma))[1 - Q(a(0)|a(\Gamma))] \geq \tfrac{1}{2} q_{\min}.$$

As $R_n$ is constant for $n$ with $l < \log|\bar{\Lambda}_n| \leq l+1$, the corresponding $Y_n$'s are actually partial sums of a sequence of $N_{n^*}^k(a(\Gamma)) \leq e^{j+1}$ such conditionally i.i.d. random variables, where $n^*$ is the largest element of $\mathcal{N}_{j,l}(k,a(\Gamma))$.

Therefore, using Lemma A.4 in the Appendix with $\mu = \mu_j = (1-\eta)\sqrt{e^{-1}\delta j^{1/2}}$, where $\eta > 0$ is an arbitrary constant, we have

$$\begin{aligned}
\mathrm{Prob}&\left\{\max_{n\in\mathcal{N}_{j,l}(k,a(\Gamma))} Y_n \ge \sqrt{e^j\delta j^{1/2}}\,\Big|\, x\Bigg(\bigcup_{i\in\bar\Lambda_n^k:\, x(\Gamma^i)=a(\Gamma)}\Gamma^i\Bigg)\right\}\\
&\le \mathrm{Prob}\left\{\max_{n\in\mathcal{N}_{j,l}(k,a(\Gamma))} Y_n \ge D\sqrt{e^{j+1}}((1-\eta)\sqrt{e^{-1}\delta j^{1/2}}+2)\right.\\
&\qquad\qquad\qquad\qquad\qquad \left.\Big|\, x\Bigg(\bigcup_{i\in\bar\Lambda_n^k:\, x(\Gamma^i)=a(\Gamma)}\Gamma^i\Bigg)\right\}\\
&\le \frac{8}{3}\exp\left[-\frac{\mu_j^2}{2(1+\mu_j/(2D\sqrt{e^{j+1}}\,))^2}\right].
\end{aligned}$$

On account of $\lim_{j\to\infty}\mu_j/(2D\sqrt{e^{j+1}}\,) = 0$, the last bound can be continued for $j > j_0$, as

$$\le \frac{8}{3}\exp\left[-\frac{(1-\eta)^2}{2e(1+\eta)}\delta j^{1/2}\right].$$

This bound also holds for the unconditional probabilities, hence we obtain from (3.6),

$$\begin{aligned}
\mathrm{Prob}\{G_j(k,a(\Gamma,0))\} &\le (\varepsilon j+2)\cdot\frac{8}{3}\exp\left[-\frac{(1-\eta)^2}{2e(1+\eta)}\delta j^{1/2}\right]\\
&\le \exp\left[-\frac{(1-\eta)^3}{2e(1+\eta)}\delta j^{1/2}\right].
\end{aligned}$$

To bound the number of all admissible $k$, $\Gamma$, $a(\Gamma,0)$ [recall the conditions $\|k\| \le 2R^{(j)}$, $r(\Gamma) \le R^{(j)}$, with $R^{(j)}$ defined in (3.5)], note that the number of possible $k$'s is bounded by

$$(4R^{(j)}+1)^d \le (4+\rho)^d\alpha^{1/2}(1+\varepsilon)^{1/2}(j+1)^{1/2},$$

and, by Lemma A.6 in the Appendix, the number of possible blocks $a(\Gamma,0)$ with $r(\Gamma) \le R^{(j)}$ is bounded by

$$(|A|^2+1)^{(2R^{(j)}+1)^d/2} < (|A|^2+1)^{(1+\rho)^d 2^{d-1}\alpha^{1/2}(1+\varepsilon)^{1/2}(j+1)^{1/2}}.$$

Combining the above bounds, we get for the probability of the union of the events $G_j(k,a(\Gamma,0))$ for all admissible $k$, $\Gamma$, $a(\Gamma,0)$ the bound

$$\begin{aligned}
\exp\Bigg[&-\frac{(1-\eta)^3}{2e(1+\eta)}\delta j^{1/2}\\
&+[\log(|A|^2+1)](1+\rho)^d 2^{d-1}\alpha^{1/2}(1+\varepsilon)^{1/2}(j+1)^{1/2} + O(\log j^{1/2})\Bigg].
\end{aligned}$$

This is summable in $j$ if we choose $\eta$, $\varepsilon$, $\rho$ sufficiently small, and $\delta/(2e) > 2^{d-1}\alpha^{1/2}\log(|A|^2+1)$, that is, if $\delta > 2^d e\alpha^{1/2}\log(|A|^2+1)$. □

Proof of Proposition 3.1. Using Lemma 3.2,

$$\begin{aligned}\left|\frac{N_n(a(\Gamma,0))}{N_n(a(\Gamma))} - Q(a(0)|a(\Gamma))\right| \\ &\le \sum_{k\colon \|k\|\le 2R_n} \left|\frac{N_n^k(a(\Gamma,0))}{N_n^k(a(\Gamma))} - Q(a(0)|a(\Gamma))\right| \cdot \frac{N_n^k(a(\Gamma))}{N_n(a(\Gamma))} \\ &< \sum_{k\colon \|k\|\le 2R_n} \sqrt{\frac{\delta\log^{1/2} N_n^k(a(\Gamma))}{N_n^k(a(\Gamma))}} \cdot \frac{N_n^k(a(\Gamma))}{N_n(a(\Gamma))}\end{aligned}$$

eventually almost surely as $n\to\infty$. By Jensen's inequality and $N_n^k(a(\Gamma)) \le N_n(a(\Gamma))$, this can be continued as

$$\le \sqrt{\frac{\delta(4R_n+1)^d\log^{1/2} N_n(a(\Gamma))}{N_n(a(\Gamma))}}.$$

By (3.4) and Lemma 3.1, we have for any $\varepsilon$, $\rho > 0$ and $n$ sufficiently large,

$$\begin{aligned}(4R_n+1)^d &\le (4\alpha^{1/(2d)}(1+\log|\bar{\Lambda}_n|)^{1/(2d)}+1)^d \\ &\le (4+\rho)^d\alpha^{1/2}(1+\varepsilon)^{1/2}\log^{1/2} N_n(a(\Gamma)),\end{aligned}$$

eventually almost surely as $n\to\infty$. This completes the proof. □

## 4. The overestimation.

Proposition 4.1. *Eventually almost surely as $n\to\infty$,*

$$\widehat{\Gamma}_{\mathrm{PIC}}(x(\Lambda_n)) \notin \{\Gamma : \Gamma \supset \Gamma_0\},$$

*whenever $r_n$ in Theorem 2.1 is equal to $R_n$ in (3.4) with*

$$\alpha < \frac{q_{\min}}{2^{3d}e}\frac{|A|-1}{|A|^2\log(|A|^2+1)}.$$

Proof. We have to prove that simultaneously for all neighborhoods $\Gamma \supset \Gamma_0$ with $r(\Gamma) \le R_n$,

$$\mathrm{PIC}_\Gamma(x(\Lambda_n)) - \mathrm{PIC}_{\Gamma_0}(x(\Lambda_n)) > 0, \tag{4.7}$$

eventually almost surely as $n\to\infty$.

The left-hand side

$$-\log \mathrm{MPL}_\Gamma(x(\Lambda_n)) + |A|^{|\Gamma|}\log|\Lambda_n| + \log \mathrm{MPL}_{\Gamma_0}(x(\Lambda_n)) - |A|^{|\Gamma_0|}\log|\Lambda_n|$$

is bounded below by

$$-\log \mathrm{MPL}_\Gamma(x(\Lambda_n)) + \log \mathrm{PL}_{\Gamma_0}(x(\Lambda_n), Q_{\Gamma_0}) + \left(1 - \frac{1}{|A|}\right)|A|^{|\Gamma|}\log|\Lambda_n|.$$

Hence, it suffices to show that simultaneously for all neighborhoods $\Gamma \supset \Gamma_0$ with $r(\Gamma) \le R_n$,

$$\log \mathrm{MPL}_\Gamma(x(\Lambda_n)) - \log \mathrm{PL}_{\Gamma_0}(x(\Lambda_n), Q_{\Gamma_0}) < \frac{|A|-1}{|A|}|A|^{|\Gamma|}\log|\Lambda_n|, \tag{4.8}$$

eventually almost surely as $n \to \infty$.

Now, for $\Gamma \supset \Gamma_0$ we have $\mathrm{PL}_{\Gamma_0}(x(\Lambda_n), Q_{\Gamma_0}) = \mathrm{PL}_\Gamma(x(\Lambda_n), Q_\Gamma)$, by the definition (2.2) of pseudo-likelihood, since $\Gamma_0$ is a Markov neighborhood. Thus the left-hand side of (4.8) equals

$$\begin{aligned}
&\log \mathrm{MPL}_\Gamma(x(\Lambda_n)) - \log \mathrm{PL}_\Gamma(x(\Lambda_n), Q_\Gamma) \\
&\quad = \sum_{a(\Gamma,0)\in x(\Lambda_n)} N_n(a(\Gamma,0)) \log \frac{N_n(a(\Gamma,0))/N_n(a(\Gamma))}{Q(a(0)|a(\Gamma))} \\
&\quad = \sum_{a(\Gamma)\in x(\Lambda_n)} N_n(a(\Gamma)) \\
&\qquad\quad \times \sum_{a(0)\,:\,a(\Gamma,0)\in x(\Lambda_n)} \frac{N_n(a(\Gamma,0))}{N_n(a(\Gamma))} \log \frac{N_n(a(\Gamma,0))/N_n(a(\Gamma))}{Q(a(0)|a(\Gamma))}.
\end{aligned}$$

To bound the last expression, we use Proposition 3.1 and Lemma A.7 in the Appendix, the latter applied with $P(a(0)) = \frac{N_n(a(\Gamma,0))}{N_n(a(\Gamma))}$, $Q(a(0)) = Q(a(0)|a(\Gamma))$. Thus we obtain the upper bound

$$\begin{aligned}
&\sum_{a(\Gamma)\in x(\Lambda_n)} N_n(a(\Gamma)) \frac{1}{q_{\min}} \sum_{a(0)\,:\,a(\Gamma,0)\in x(\Lambda_n)} \left[\frac{N_n(a(\Gamma,0))}{N_n(a(\Gamma))} - Q(a(0)|a(\Gamma))\right]^2 \\
&\quad < \sum_{a(\Gamma)\in x(\Lambda_n)} N_n(a(\Gamma)) \frac{1}{q_{\min}} |A| \frac{\kappa \log N_n(a(\Gamma))}{N_n(a(\Gamma))} \le \frac{\kappa|A|}{q_{\min}}|A|^{|\Gamma|}\log|\bar{\Lambda}_n|,
\end{aligned}$$

eventually almost surely as $n \to \infty$, simultaneously for all neighborhoods $\Gamma \supset \Gamma_0$ with $r(\Gamma) \le R_n$.

Hence, since $|\Lambda_n|/|\bar{\Lambda}_n| \to 1$, the assertion (4.8) holds whenever

$$\frac{\kappa|A|}{q_{\min}} < \frac{|A|-1}{|A|},$$

which is equivalent to the bound on $\alpha$ in Proposition 4.1. □

## 5. The underestimation.

Proposition 5.1. *Eventually almost surely as $n\to\infty$,*

$$\widehat{\Gamma}_{\mathrm{PIC}}(x(\Lambda_n))\in\{\Gamma:\Gamma\supseteq\Gamma_0\},$$

*if $r_n$ in Theorem 2.1 is chosen as in Proposition 4.1.*

Proposition 5.1 will be proved using the lemmas below. Let us denote

$$\Psi_0=\left(\bigcup_{i\in\Gamma_0}\Gamma_0^i\right)\Big\backslash(\Gamma_0\cup\{0\}).$$

Lemma 5.1. *The assertion of Proposition 3.1 holds also with $\Gamma$ replaced by $\Gamma\cup\Psi_0$, where $\Gamma$ is any (not necessarily Markov) neighborhood.*

Proof. As Proposition 3.1 was a consequence of Lemma 3.2, we have to check that the proof of that lemma works when the Markov neighborhood $\Gamma$ is replaced by $\Gamma\cup\Psi_0$, where $\Gamma$ is any neighborhood. To this end, it suffices to show that conditional on the values of all random variables in the $(\Gamma\cup\Psi_0)$-neighborhoods of the sites $i\in\bar{\Lambda}_n^k$, those $X(i)$, $i\in\bar{\Lambda}_n^k$, are conditionally i.i.d. for which the same block $a(\Gamma\cup\Psi_0)$ appears in the $(\Gamma\cup\Psi_0)$-neighborhood of $i$. This follows from Lemma A.1 in the Appendix, with $\Delta=\Gamma_0\cup\{0\}$ and $\Psi=\Psi_0$. □

Lemma 5.2. *Simultaneously for all neighborhoods $\Gamma\not\supseteq\Gamma_0$ with $r(\Gamma)\le R_n$,*

$$\mathrm{PIC}_{\Gamma\cup\Psi_0}(x(\Lambda_n))>\mathrm{PIC}_{(\Gamma\cap\Gamma_0)\cup\Psi_0}(x(\Lambda_n)),$$

*eventually almost surely as $n\to\infty$.*

Proof. The claimed inequality is analogous to (4.7) in the proof of Proposition 4.1, the role of $\Gamma\supset\Gamma_0$ there played by $\Gamma\cup\Psi_0\supset(\Gamma\cap\Gamma_0)\cup\Psi_0$. Its proof is the same as that of (4.7), using Lemma 5.1 instead of Proposition 3.1. Indeed, the basic neighborhood property of $\Gamma_0$ was used in that proof only to show that $\mathrm{PL}_{\Gamma_0}(x(\Lambda_n),Q_{\Gamma_0})=\mathrm{PL}_{\Gamma}(x(\Lambda_n),Q_{\Gamma})$. The analogue of this identity, namely

$$\mathrm{PL}_{(\Gamma\cap\Gamma_0)\cup\Psi_0}(x(\Lambda_n),Q_{(\Gamma\cap\Gamma_0)\cup\Psi_0})=\mathrm{PL}_{\Gamma\cup\Psi_0}(x(\Lambda_n),Q_{\Gamma\cup\Psi_0}),$$

follows from Lemma A.1 in the Appendix with $\Delta=\Gamma_0\cup\{0\}$ and $\Psi=\Psi_0$. □

For the next lemma, we introduce some further notation.

The set of all probability distributions on $A^{\mathbb{Z}^d}$, equipped with the weak topology, is a compact Polish space; let $d$ denote a metric that metrizes it. Let $\mathcal{Q}^G$ denote the (compact) set of Gibbs distributions with the one-point specification $Q_{\Gamma_0}$.

For a sample $x(\Lambda_n)$, define the empirical distribution on $A^{\mathbb{Z}^d}$ by

$$R_{x,n} = \frac{1}{|\bar{\Lambda}_n|} \sum_{i \in \bar{\Lambda}_n} \delta_{x_n^i},$$

where $x_n \in A^{\mathbb{Z}^d}$ is the extension of the sample $x(\Lambda_n)$ to the whole lattice with $x_n(j)$ equal to a constant $a \in A$ for $j \in \mathbb{Z}^d \backslash \Lambda_n$, and $x_n^i$ denotes the translate of $x_n$ when 0 is translated to $i$ and $\delta_x$ is the Dirac mass at $x \in A^{\mathbb{Z}^d}$.

LEMMA 5.3. *With probability* 1, $d(R_{x,n}, \mathcal{Q}^G) \to 0$.

PROOF. Fix a realization $x(\mathbb{Z}^d)$ for which Proposition 3.1 holds.

It suffices to show that for any subsequence $n_k$ such that $R_{x,n_k}$ converges, its limit $R_{x,0}$ belongs to $\mathcal{Q}^G$.

Let $\Gamma'$ be any neighborhood. For $n$ sufficiently large, the $(\Gamma' \cup \{0\})$-marginal of $R_{x,n}$ is equal to

$$\left\{ \frac{N_n(a(\Gamma', 0))}{|\bar{\Lambda}_n|}, a(\Gamma', 0) \in A^{\Gamma' \cup \{0\}} \right\},$$

hence $R_{x,n_k} \to R_{x,0}$ implies

$$\frac{N_{n_k}(a(\Gamma', 0))}{|\bar{\Lambda}_{n_k}|} \longrightarrow R_{x,0}(a(\Gamma', 0)) \tag{5.9}$$

for all $a(\Gamma', 0) \in A^{\Gamma' \cup \{0\}}$. This and summation for $a(0) \in A$ imply

$$\frac{N_{n_k}(a(\Gamma', 0))}{N_{n_k}(a(\Gamma'))} \longrightarrow R_{x,0}(a(0)|a(\Gamma')).$$

As Proposition 3.1 holds for the realization $x(\mathbb{Z}^d)$, it follows that if $\Gamma'$ is a Markov neighborhood, then

$$R_{x,0}(a(0)|a(\Gamma')) = Q(a(0)|a(\Gamma')) = Q_{\Gamma_0}(a(0)|a(\Gamma_0)).$$

For any finite region $\Delta \supset \Gamma_0$ with $0 \notin \Delta$, the last equation for a neighborhood $\Gamma' \supset \Delta$ implies that

$$R_{x,0}(a(0)|a(\Delta)) = Q_{\Gamma_0}(a(0)|a(\Gamma_0)) \qquad \text{if } \Delta \supset \Gamma_0, 0 \notin \Delta.$$

To prove $R_{x,0} \in \mathcal{Q}^G$ it remains to show that, in addition, $R_{x,0}(a(i)|a(\Delta^i)) = Q_{\Gamma_0}(a(i)|a(\Gamma_0^i))$. Actually, we show that $R_{x,0}$ is translation invariant. Indeed, given a finite region $\Delta \subset \mathbb{Z}^d$ and its translate $\Delta^i$, take a neighborhood $\Gamma'$

with $\Delta \cup \Delta^i \subseteq \Gamma' \cup \{0\}$, and consider the sum of the counts $N_n(a(\Gamma', 0))$ for all blocks $a(\Gamma', 0) = \{a(j) : j \in \Gamma' \cup \{0\}\}$ with $\{a(j) : j \in \Delta\}$ equal to a fixed $|\Delta|$-tuple and the similar sum with $\{a(j) : j \in \Delta^i\}$ equal to the same $|\Delta|$-tuple. If $\|i\| < \log^{1/(2d)} |\bar{\Lambda}_n|$, the difference of these sums is at most $|\Lambda_n| - |\bar{\Lambda}_n|$, hence the translation invariance of $R_{x,0}$ follows by (5.9). □

LEMMA 5.4. *Uniformly for all neighborhoods $\Gamma$ not containing $\Gamma_0$,*

$$-\log \mathrm{MPL}_{(\Gamma\cap\Gamma_0)\cup\Psi_0}(x(\Lambda_n)) > -\log \mathrm{MPL}_{\Gamma_0}(x(\Lambda_n)) + c|\bar{\Lambda}_n|,$$

*eventually almost surely as $n \to \infty$, where $c > 0$ is a constant.*

PROOF. Given a realization $x \in A^{\mathbb{Z}^d}$ with the property in Lemma 5.3, there exists a sequence $Q_{x,n}$ in $\mathcal{Q}^G$ with

$$d(R_{x,n}, Q_{x,n}) \to 0,$$

and consequently

$$\frac{N_n(a(\Delta))}{|\bar{\Lambda}_n|} - Q_{x,n}(a(\Delta)) \to 0 \tag{5.10}$$

for each finite region $\Delta \subset \mathbb{Z}^d$ and $a(\Delta) \in A^\Delta$.

Next, let $\Gamma$ be a neighborhood with $\Gamma \not\supseteq \Gamma_0$. By (2.3),

$$\begin{aligned}
-&\frac{1}{|\bar{\Lambda}_n|}\log \mathrm{MPL}_{(\Gamma\cap\Gamma_0)\cup\Psi_0}(x(\Lambda_n)) \\
&= -\frac{1}{|\bar{\Lambda}_n|}\sum_{a((\Gamma\cap\Gamma_0)\cup\Psi_0,0)\in x(\Lambda_n)} N_n(a((\Gamma\cap\Gamma_0)\cup\Psi_0, 0)) \\
&\qquad\qquad \times \log\frac{N_n(a((\Gamma\cap\Gamma_0)\cup\Psi_0,0))}{N_n(a((\Gamma\cap\Gamma_0)\cup\Psi_0))}.
\end{aligned}$$

Applying (5.10) to $\Delta = (\Gamma\cap\Gamma_0)\cup\Psi_0\cup\{0\}$, it follows that the last expression is arbitrarily close to

$$\begin{aligned}
-&\sum_{a((\Gamma\cap\Gamma_0)\cup\Psi_0\cup\{0\})} Q_{x,n}(a((\Gamma\cap\Gamma_0)\cup\Psi_0,0))\log Q_{x,n}(a(0)|a((\Gamma\cap\Gamma_0)\cup\Psi_0)) \\
&= H_{Q_{x,n}}(X(0)|X((\Gamma\cap\Gamma_0)\cup\Psi_0))
\end{aligned}$$

if $n$ is sufficiently large, where $H_{Q_{x,n}}(\cdot|\cdot)$ denotes conditional entropy, when the underlying distribution is $Q_{x,n}$. Similarly, $-(1/|\bar{\Lambda}_n|)\log \mathrm{MPL}_{\Gamma_0}(x(\Lambda_n))$ is arbitrarily close to $H_{Q_{x,n}}(X(0)|X(\Gamma_0))$, which equals $H_{Q_{x,n}}(X(0)|X(\Gamma_0 \cup \Psi_0))$ since $\Gamma_0$ is a Markov neighborhood.

It is known that $H_{Q'}(X(0)|X((\Gamma\cap\Gamma_0)\cup\Psi_0)) \geq H_{Q'}(X(0)|X(\Gamma_0\cup\Psi_0))$ for any distribution $Q'$. The proof of the lemma will be complete if we show

that, in addition, there exists a constant $\xi > 0$ (depending on $\Gamma \cap \Gamma_0$) such that for every Gibbs distribution $Q^G \in \mathcal{Q}^G$

$$H_{Q^G}(X(0)|X((\Gamma \cap \Gamma_0) \cup \Psi_0)) - H_{Q^G}(X(0)|X(\Gamma_0 \cup \Psi_0)) > \xi.$$

The indirect assumption that the left-hand side goes to 0 for some sequence of Gibbs distributions in $\mathcal{Q}^G$ implies, using the compactness of $\mathcal{Q}^G$, that

$$H_{Q_0^G}(X(0)|X((\Gamma \cap \Gamma_0) \cup \Psi_0)) = H_{Q_0^G}(X(0)|X(\Gamma_0 \cup \Psi_0)),$$

for the limit $Q_0^G \in \mathcal{Q}^G$ of a convergent subsequence. This equality implies

$$Q_0^G(a(0)|a((\Gamma \cap \Gamma_0) \cup \Psi_0)) = Q_0^G(a(0)|a(\Gamma_0 \cup \Psi_0))$$

for all $a(0) \in A, a(\Gamma_0 \cup \Psi_0) \in A^{\Gamma_0 \cup \Psi_0}$. By Lemma A.1 in the Appendix, these conditional probabilities are uniquely determined by the one-point specification $Q_{\Gamma_0}$, and the last equality implies

$$Q(a(i)|a((\Gamma \cap \Gamma_0)^i \cup \Psi_0^i)) = Q(a(i)|a(\Gamma_0^i \cup \Psi_0^i)) = Q_{\Gamma_0}(a(i)|a(\Gamma_0^i)).$$

According to Lemma A.2 in the Appendix, this would imply $(\Gamma \cap \Gamma_0) \cup \Psi_0$ is a Markov neighborhood also, which is a contradiction, as $(\Gamma \cap \Gamma_0) \cup \Psi_0 \not\supseteq \Gamma_0$.

This completes the proof of the lemma because there is only a finite number of possible intersections $\Gamma \cap \Gamma_0$. □

Proof of Proposition 5.1. We have to show that

$$\mathrm{PIC}_\Gamma(x(\Lambda_n)) > \mathrm{PIC}_{\Gamma_0}(x(\Lambda_n)), \tag{5.11}$$

eventually almost surely as $n \to \infty$, for all neighborhoods $\Gamma$ with $r(\Gamma) \leq R_n$ that do not contain $\Gamma_0$.

Note that $\Gamma_1 \supseteq \Gamma_2$ implies $\mathrm{MPL}_{\Gamma_1}(x(\Lambda_n)) \geq \mathrm{MPL}_{\Gamma_2}(x(\Lambda_n))$, since $\mathrm{MPL}_\Gamma(x(\Lambda_n))$ is the maximizer in $Q'_\Gamma$ of $\mathrm{PL}_\Gamma(x(\Lambda_n), Q'_\Gamma)$; see (2.2). Hence

$$-\log \mathrm{MPL}_\Gamma(x(\Lambda_n)) \geq -\log \mathrm{MPL}_{\Gamma \cup \Psi_0}(x(\Lambda_n))$$

for any neighborhood $\Gamma$.

Thus

$$\begin{aligned}\mathrm{PIC}_\Gamma(x(\Lambda_n)) &= -\log \mathrm{MPL}_\Gamma(x(\Lambda_n)) + |A|^{|\Gamma|} \log |\Lambda_n| \\ &\geq \mathrm{PIC}_{\Gamma \cup \Psi_0}(x(\Lambda_n)) - (|A|^{|\Gamma \cup \Psi_0|} - |A|^{|\Gamma|}) \log |\Lambda_n|.\end{aligned}$$

Using Lemma 5.2 and the obvious bound $|\Gamma \cup \Psi_0| \leq |\Gamma| + |\Psi_0|$, it follows that, eventually almost surely as $n \to \infty$ for all $\Gamma \not\supseteq \Gamma_0$ with $r(\Gamma) \leq R_n$,

$$\mathrm{PIC}_\Gamma(x(\Lambda_n)) > \mathrm{PIC}_{(\Gamma \cap \Gamma_0) \cup \Psi_0}(x(\Lambda_n)) - |A|^{|\Gamma|}(|A|^{|\Psi_0|} - 1) \log |\Lambda_n|.$$

Here, by Lemma 5.4,

$$\begin{aligned}&\mathrm{PIC}_{(\Gamma \cap \Gamma_0) \cup \Psi_0}(x(\Lambda_n)) \\ &\qquad > -\log \mathrm{MPL}_{(\Gamma \cap \Gamma_0) \cup \Psi_0}(x(\Lambda_n)) > -\log \mathrm{MPL}_{\Gamma_0}(x(\Lambda_n)) + c|\bar{\Lambda}_n|,\end{aligned}$$

eventually almost surely as $n \to \infty$ for all $\Gamma$ as above. This completes the proof, since the conditions $r(\Gamma) \leq R_n$ and $|\Lambda_n|/|\bar{\Lambda}_n| \to 1$ imply $|A|^{|\Gamma|} \log |\Lambda_n| = o(|\bar{\Lambda}_n|)$. □

**6. Discussion.** A modification of the Bayesian Information Criterion (BIC) called PIC has been introduced for estimating the basic neighborhood of a Markov random field on $\mathbb{Z}^d$, with finite alphabet $A$. In this criterion, the maximum pseudo-likelihood is used instead of the maximum likelihood, with penalty term $|A|^{|\Gamma|} \log |\Lambda_n|$ for a candidate neighborhood $\Gamma$, where $\Lambda_n$ is the sample region. The minimizer of PIC over candidate neighborhoods, with radius allowed to grow as $o(\log^{1/(2d)} |\Lambda_n|)$, has been proved to equal the basic neighborhood eventually almost surely, not requiring any prior bound on the size of the latter. This result is unaffected by phase transition and even by nonstationarity of the joint distribution. The same result holds if the penalty term is multiplied by any $c > 0$; the no underestimation part (Proposition 5.1) holds also if $\log |\Lambda_n|$ in the penalty term is replaced by any function of the sample size $|\Lambda_n|$ that goes to infinity as $o(|\Lambda_n|)$.

PIC estimation of the basic neighborhood of a Markov random field is to a certain extent similar to BIC estimation of the order of a Markov chain, and of the context tree of a tree source, also called a variable-length Markov chain. For context tree estimation via another method see [5, 22], and via BIC, see [9]. There are, however, also substantial differences. The martingale techniques in [7, 8] do not appear to carry over to Markov random fields, and the lack of an analogue of the Krichevsky–Trofimov distribution used in these references is another obstacle. We also note that the "large" boundaries of multidimensional sample regions cause side effects not present in the one-dimensional case; to overcome those, we have defined the pseudo-likelihood function based on a window $\bar{\Lambda}_n$ slightly smaller than the whole sample region $\Lambda_n$.

For Markov order and context tree estimation via BIC, consistency has been proved by Csiszár and Shields [8] admitting, for sample size $n$, all $k \leq n$ as candidate orders (see also [7]), respectively by Csiszár and Talata [9] admitting trees of depth $o(\log n)$ as candidate context trees. In our main result Theorem 2.1, the PIC estimator of the basic neighborhood is defined admitting candidate neighborhoods of radius $o(\log^{1/(2d)} |\Lambda_n|)$, thus of size $o(\log^{1/2} |\Lambda_n|)$. The mentioned one-dimensional results suggest that this bound on the radius might be relaxed to $o(\log^{1/d} |\Lambda_n|)$, or perhaps dropped completely. This question remains open, even for the case $d = 1$. A positive answer apparently depends on the possibility of strengthening our typicality result Proposition 3.1 to similar strength as the conditional typicality results for Markov chains in [7].

More important than a possible mathematical sharpening of Theorem 2.1, as above, would be to find an algorithm to determine the PIC estimator without actually computing and comparing the PIC values of all candidate neighborhoods. The analogous problem for BIC context tree estimation has been solved: Csiszár and Talata [9] showed that this BIC estimator can be computed in linear time via an analogue of the "context tree maximizing algorithm" of Willems, Shtarkov and Tjalkens [23, 24]. Unfortunately, a similar algorithm for the present problem appears elusive, and it remains open whether our estimator can be computed in a "clever" way.

Finally, we emphasize that the goal of this paper was to provide a consistent estimator of the basic neighborhood of a Markov random field. Of course, consistency is only one of the desirable properties of an estimator. To assess the practical performance of this estimator requires further research, such as studying finite sample size properties, robustness against noisy observations and computability with acceptable complexity.

*Note added in proof.* Just before completing the galley proofs, we learned that model selection for Markov random fields had been addressed before, by Ji and Seymour [18]. They used a criterion almost identical to PIC here and, in a somewhat different setting, proved weak consistency under the assumption that the number of candidate model classes is finite.

## APPENDIX

First we indicate how the well-known facts stated in Lemma 2.1 can be formally derived from results in [13], using the concepts defined there.

Proof of Lemma 2.1. By Theorem 1.33 the positive one-point specification uniquely determines the specification, which is positive and local on account of the locality of the one-point specification. By Theorem 2.30 this positive local specification determines a unique "gas" potential (if an element of $A$ is distinguished as the zero element). Due to Corollary 2.32, this is a nearest-neighbor potential for a graph with vertex set $\mathbb{Z}^d$ defined there, and $\Gamma_0^i$ is the same as $B(i)\backslash\{i\}$ in that corollary. □

The following lemma is a consequence of the global Markov property.

Lemma A.1. *Let $\Delta \subset \mathbb{Z}^d$ be a finite region with $0 \in \Delta$, and $\Psi = (\bigcup_{j\in\Delta} \Gamma_0^j) \setminus \Delta$. Then for any neighborhood $\Gamma$, the conditional probabilities $Q(a(i)|a(\Gamma^i \cup \Psi^i))$ and $Q(a(i)|a((\Gamma^i \cap \Delta^i) \cup \Psi^i))$ are equal and translation invariant.*

Proof. Since $\Delta$ and $\Psi$ are disjoint, we have

$$\begin{aligned} Q(a(i)|a(\Gamma^i \cup \Psi^i)) &= Q(a(i)|a((\Gamma \cap \Delta)^i \cup (\Psi \cup (\Gamma\backslash\Delta))^i)) \\ &= \frac{Q(a(\{i\} \cup (\Gamma \cap \Delta)^i)|a((\Psi \cup (\Gamma\backslash\Delta))^i))}{Q(a((\Gamma \cap \Delta)^i)|a((\Psi \cup (\Gamma\backslash\Delta))^i))}, \end{aligned}$$

and similarly

$$Q(a(i)|a((\Gamma^i \cap \Delta^i) \cup \Psi^i)) = \frac{Q(a(\{i\} \cup (\Gamma \cap \Delta)^i)|a(\Psi^i))}{Q(a((\Gamma \cap \Delta)^i)|a(\Psi^i))}.$$

By the global Markov property (see Lemma 2.1), both the numerators and denominators of these two quotients are equal and translation invariant. □

The lemma below follows from the definition of Markov neighborhood.

LEMMA A.2. *For a Markov random field with basic neighborhood $\Gamma_0$, if a neighborhood $\Gamma$ satisfies*

$$Q(a(i)|a(\Gamma^i)) = Q_{\Gamma_0}(a(i)|a(\Gamma_0^i))$$

*for all $i \in \mathbb{Z}^d$, then $\Gamma$ is a Markov neighborhood.*

PROOF. We have to show that for any $\Delta \supset \Gamma$

$$Q(a(i)|a(\Delta^i)) = Q(a(i)|a(\Gamma^i)). \tag{A.1}$$

Since $\Gamma_0$ is a Markov neighborhood, the condition of the lemma implies

$$Q(a(i)|a(\Gamma^i)) = Q(a(i)|a(\Gamma_0^i)) = Q(a(i)|a((\Gamma_0 \cup \Delta)^i)).$$

Hence (A.1) follows, because $\Gamma \subseteq \Delta \subseteq \Gamma_0 \cup \Delta$. □

Next we state two simple probability bounds.

LEMMA A.3. *Let $Z_1, Z_2, \ldots$ be $\{0,1\}$-valued random variables such that*

$$\text{Prob}\{Z_j = 1|Z_1, \ldots Z_{j-1}\} \geq p_* > 0, \qquad j \geq 1,$$

*with probability 1. Then for any $0 < \nu < 1$*

$$\text{Prob}\left\{\frac{1}{m}\sum_{j=1}^{m} Z_j < \nu p_*\right\} \leq e^{-m(p_*/4)(1-\nu)^2}.$$

PROOF. This is a direct consequence of Lemmas 2 and 3 in the Appendix of [7]. □

LEMMA A.4. *Let $Z_1, Z_2, \ldots, Z_n$ be i.i.d. random variables with expectation 0 and variance $D^2$. Then the partial sums*

$$S_k = Z_1 + Z_2 + \cdots + Z_k$$

*satisfy*

$$\text{Prob}\left\{\max_{1 \leq k \leq n} S_k \geq D\sqrt{n}(\mu + 2)\right\} \leq \tfrac{4}{3}\text{Prob}\{S_n \geq D\sqrt{n}\mu\};$$

*moreover if the random variables are bounded, $|Z_i| \le K$, then*

$$\mathrm{Prob}\{S_n \ge D\sqrt{n}\mu\} \le 2\exp\left[-\frac{\mu^2}{2(1+\mu K(2D\sqrt{n}\,))^2}\right],$$

*where $\mu < D\sqrt{n}/K$.*

Proof. See, for example, Lemma VI.9.1 and Theorem VI.4.1 in [20]. □

The following three lemmas are of a technical nature.

Lemma A.5. *For disjoint finite regions $\Phi \subset \mathbb{Z}^d$ and $\Delta \subset \mathbb{Z}^d$, we have*

$$Q(a(\Delta)|a(\Phi)) \ge q_{\min}^{|\Delta|}.$$

Proof. By induction on $|\Delta|$.
For $\Delta = \{i\}$, $\Xi = \Gamma_0^i \setminus \Phi$, we have

$$\begin{aligned} Q(a(i)|a(\Phi)) &= \sum_{a(\Xi)\in A^\Xi} Q(a(i)|a(\Phi\cup\Xi))Q(a(\Xi)|a(\Phi)) \\ &= \sum_{a(\Xi)\in A^\Xi} Q(a(i)|a(\Gamma_0^i))Q(a(\Xi)|a(\Phi)) \ge q_{\min}. \end{aligned}$$

Supposing $Q(a(\Delta)|a(\Phi)) \ge q_{\min}^{|\Delta|}$ holds for some $\Delta$, we have for $\{i\}\cup\Delta$, with $\Xi = \Gamma_0^i \setminus (\Phi\cup\Delta)$,

$$\begin{aligned} Q(a(\{i\}\cup\Delta)|a(\Phi)) &= \sum_{a(\Xi)\in A^\Xi} Q(a(\{i\}\cup\Delta\cup\Xi)|a(\Phi)) \\ &= \sum_{a(\Xi)\in A^\Xi} Q(a(i)|a(\Delta\cup\Xi\cup\Phi))Q(a(\Delta\cup\Xi)|a(\Phi)). \end{aligned}$$

Since $Q(a(i)|a(\Delta\cup\Xi\cup\Phi)) = Q(a(i)|a(\Gamma_0^i)) \ge q_{\min}$, we can continue as

$$\ge q_{\min} Q(a(\Delta)|a(\Phi)) \ge q_{\min}^{|\Delta|+1}. \qquad \square$$

Lemma A.6. *The number of all possible blocks appearing in a site and its neighborhood with radius not exceeding $R$ can be upper bounded as*

$$|\{a(\Gamma,0) \in A^{\Gamma\cup\{0\}} : r(\Gamma) \le R\}| \le (|A|^2+1)^{(2R+1)^d/2}.$$

Proof. The number of the neighborhoods with cardinality $m \ge 1$ and radius $r(\Gamma) \le R$ is

$$\binom{((2R+1)^d-1)/2}{m},$$

because the neighborhoods are symmetric. Hence, the number in the proposition is

$$|A| + |A| \cdot \sum_{m=1}^{((2R+1)^d-1)/2} \binom{((2R+1)^d-1)/2}{m} |A|^{2m}$$
$$= |A| \sum_{m=0}^{((2R+1)^d-1)/2} \binom{((2R+1)^d-1)/2}{m} (|A|^2)^m 1^{((2R+1)^d-1)/2-m}.$$

Now, using the binomial theorem, the assertion follows. □

Lemma A.7. *Let $P$ and $Q$ be probability distributions on $A$ such that*

$$\max_{a\in A} |P(a) - Q(a)| \le \frac{\min_{a\in A} Q(a)}{2}.$$

*Then*

$$\sum_{a\in A} P(a) \log \frac{P(a)}{Q(a)} \le \frac{1}{\min_{a\in A} Q(a)} \sum_{a\in A} (P(a) - Q(a))^2.$$

Proof. This follows from Lemma 4 in the Appendix of [7]. □

## REFERENCES

[1] Akaike, H. (1972). Information theory and an extension of the maximum likelihood principle. In *Proc. Second International Symposium on Information Theory. Supplement to Problems of Control and Information Theory* (B. N. Petrov and F. Csáki, eds.) 267–281. Akadémiai Kiadó, Budapest. MR0483125

[2] Azencott, R. (1987). Image analysis and Markov fields. In *Proc. First International Conference on Industrial and Applied Mathematics, Paris* (J. McKenna and R. Temen, eds.) 53–61. SIAM, Philadelphia. MR0976851

[3] Besag, J. (1974). Spatial interaction and the statistical analysis of lattice systems (with discussion). *J. Roy. Statist. Soc. Ser. B* **36** 192–236. MR0373208

[4] Besag, J. (1975). Statistical analysis of non-lattice data. *The Statistician* **24** 179–195.

[5] Bühlmann, P. and Wyner, A. J. (1999). Variable length Markov chains. *Ann. Statist.* **27** 480–513. MR1714720

[6] Comets, F. (1992). On consistency of a class of estimators for exponential families of Markov random fields on the lattice. *Ann. Statist.* **20** 455–468. MR1150354

[7] Csiszár, I. (2002). Large-scale typicality of Markov sample paths and consistency of MDL order estimators. *IEEE Trans. Inform. Theory* **48** 1616–1628. MR1909476

[8] Csiszár, I. and Shields, P. C. (2000). The consistency of the BIC Markov order estimator. *Ann. Statist.* **28** 1601–1619. MR1835033

[9] Csiszár, I. and Talata, Zs. (2006). Context tree estimation for not necessarily finite memory processes, via BIC and MDL. *IEEE Trans. Inform. Theory* **52** 1007–1016.

[10] Dobrushin, R. L. (1968). The description of a random field by means of conditional probabilities and conditions for its regularity. *Theory Probab. Appl.* **13** 197–224. MR0231434

[11] Finesso, L. (1992). Estimation of the order of a finite Markov chain. In *Recent Advances in Mathematical Theory of Systems, Control, Networks and Signal Processing* **1** (H. Kimura and S. Kodama, eds.) 643–645. Mita, Tokyo. MR1197985

[12] Geman, S. and Graffigne, C. (1987). Markov random field image models and their applications to computer vision. In *Proc. International Congress of Mathematicians* **2** (A. M. Gleason, ed.) 1496–1517. Amer. Math. Soc., Providence, RI. MR0934354

[13] Georgii, H. O. (1988). *Gibbs Measures and Phase Transitions.* de Gruyter, Berlin. MR0956646

[14] Gidas, B. (1988). Consistency of maximum likelihood and pseudolikelihood estimators for Gibbs distributions. In *Stochastic Differential Systems, Stochastic Control Theory and Applications* (W. Fleming and P.-L. Lions, eds.) 129–145. Springer, New York. MR0934721

[15] Gidas, B. (1993). Parameter estimation for Gibbs distributions from fully observed data. In *Markov Random Fields*: *Theory and Application* (R. Chellappa and A. Jain, eds.) 471–498. Academic Press, Boston. MR1214376

[16] Hannan, E. J. and Quinn, B. G. (1979). The determination of the order of an autoregression. *J. Roy. Statist. Soc. Ser. B* **41** 190–195. MR0547244

[17] Haughton, D. (1988). On the choice of a model to fit data from an exponential family. *Ann. Statist.* **16** 342–355. MR0924875

[18] Ji, C. and Seymour, L. (1996). A consistent model selection procedure for Markov random fields based on penalized pseudolikelihood. *Ann. Appl. Probab.* **6** 423–443. MR1398052

[19] Pickard, D. K. (1987). Inference for discrete Markov fields: The simplest non-trivial case. *J. Amer. Statist. Assoc.* **82** 90–96. MR0883337

[20] Rényi, A. (1970). *Probability Theory.* North-Holland, Amsterdam. MR0315747

[21] Schwarz, G. (1978). Estimating the dimension of a model. *Ann. Statist.* **6** 461–464. MR0468014

[22] Weinberger, M. J., Rissanen, J. and Feder, M. (1995). A universal finite memory source. *IEEE Trans. Inform. Theory* **41** 643–652.

[23] Willems, F. M. J., Shtarkov, Y. M. and Tjalkens, T. J. (1993). The context-tree weighting method: Basic properties. Technical report, Dept. Electrical Engineering, Eindhoven Univ.

[24] Willems, F. M. J., Shtarkov, Y. M. and Tjalkens, T. J. (2000). Context-tree maximizing. In *Proc. 2000 Conf. Information Sciences and Systems* TP6-7–TP6-12. Princeton, NJ.

A. Rényi Institute of Mathematics
Hungarian Academy of Sciences
POB 127, H-1364 Budapest
Hungary
E-mail: csiszar@renyi.hu
zstalata@renyi.hu
URL: www.renyi.hu/~csiszar
www.renyi.hu/~zstalata